\documentclass[12pt]{article}
\usepackage{latexsym}
\usepackage{amsmath}
\usepackage{amssymb}
\usepackage{amscd}
\usepackage{graphicx}

\newtheorem{Th}{Theorem}[section]
\newtheorem{Cor}{Corollary}[section]
\newtheorem{Prop}{Proposition}[section]
\newtheorem{Lem}{Lemma}[section]
\newtheorem{Def}{Definition}[section]
\newtheorem{Rem}{Remark}[section]
\newtheorem{Ex}{Example}[section]

\newcommand{\bet}{\begin{Th}}
\newcommand{\ent}{\stepcounter{Cor}
   \stepcounter{Prop}\stepcounter{Lem}\stepcounter{Def}
   \stepcounter{Rem}\stepcounter{Ex}\end{Th}}

\newcommand{\bec}{\begin{Cor}}
\newcommand{\enc}{\stepcounter{Th}
   \stepcounter{Prop}\stepcounter{Lem}\stepcounter{Def}
   \stepcounter{Rem}\stepcounter{Ex}\end{Cor}}
\newcommand{\bep}{\begin{Prop}}
\newcommand{\enp}{\stepcounter{Th}
   \stepcounter{Cor}\stepcounter{Lem}\stepcounter{Def}
   \stepcounter{Rem}\stepcounter{Ex}\end{Prop}}
\newcommand{\bel}{\begin{Lem}}
\newcommand{\enl}{\stepcounter{Th}
   \stepcounter{Cor}\stepcounter{Prop}\stepcounter{Def}
   \stepcounter{Rem}\stepcounter{Ex}\end{Lem}}
\newcommand{\bef}{\begin{Def}}
\newcommand{\enf}{\stepcounter{Th}
   \stepcounter{Cor}\stepcounter{Prop}\stepcounter{Lem}
   \stepcounter{Rem}\stepcounter{Ex}\end{Def}}
\newcommand{\ber}{\begin{Rem} }
\newcommand{\enr}{
   \stepcounter{Th}\stepcounter{Cor}\stepcounter{Prop}
   \stepcounter{Lem}\stepcounter{Def}\stepcounter{Ex}\end{Rem}}
\newcommand{\bee}{\begin{Ex}}
\newcommand{\ene}{
   \stepcounter{Th}\stepcounter{Cor}\stepcounter{Prop}
   \stepcounter{Lem}\stepcounter{Def}\stepcounter{Rem}\end{Ex}}
\newcommand{\Proof}{\noindent{\it Proof\,}:\ }

\newcommand{\R}{\mathbf{R}}

\newcommand{\NN}{\mathbf{N}}

\newcommand{\pa}{\partial}

\newcommand{\xxx}{\boldsymbol{x}}
\newcommand{\nn}{\boldsymbol{n}}
\newcommand{\ee}{\boldsymbol{e}}

\newcommand{\QED}{\hfill $\Box$ \par}
\newcommand{\id}{{\mbox {\rm id}}}

\newcommand{\rank}{{\mbox {\rm rank}}}

\newcommand{\mathsym}[1]{{}}

\begin{document}

\title{
Singularities of flat extensions 
\\
from generic surfaces with boundaries 
}
\author{Goo ISHIKAWA
\\
Department of Mathematics, Hokkaido University
\\
Sapporo 060-0810, JAPAN. 
\\
E-mail : ishikawa@math.sci.hokudai.ac.jp
}

\date{
} 

\maketitle


\begin{abstract}
We solve the problem on flat extensions 
of a generic surface with boundary in Euclidean $3$-space, relating it to 
the singularity theory of the envelope generated by the boundary. 
We give related results on Legendre surfaces with boundaries via projective duality 
and observe the duality on boundary singularities. 
Moreover we give formulae related to
remote singularities of the boundary-envelope. 

\end{abstract}

\section{Introduction.}

\

We mean by {\it the flat extension problem} the problem on 
the existence, uniqueness and singularities of extensions 
of a surface across its boundary by flat surfaces in Euclidean $3$-space $\R^3$: 
\vspace{0.2truecm}

\noindent 
{\bf Problem:} 
Let $(S, \gamma)$ be a $C^\infty$ surface with boundary $\gamma$ in $\R^3$. 
Find a {\it $C^1$ extension} $\widetilde{S}$ of $S$ such that 
$\widetilde{S} \setminus {\mbox{\rm{Int}}}S$ is $C^\infty$ and the Gaussian curvature
$$
K \vert_{ \widetilde{S} \setminus {\mbox{\footnotesize\rm{Int}}}S } \ \equiv \ 0. 
$$

We call $\widetilde{S}$ a flat $C^1$ extension of $S$. Then the surface $(S, \gamma)$ with boundary 
is extended by a flat surface $(S', \gamma) = (\widetilde{S} \setminus {\mbox{\footnotesize\rm{Int}}}S, \gamma)$ 
with boundary. 
Recall that a surface $S'$ in $\R^3$ is called {\it flat} if it is locally isometric to the plane and 
the condition is equivalent to that $K\vert_{S'} = 0$ (\cite{Spivak}). 
Note that, in general, for a hypersurface $y = f(x_1, \dots, x_n)$ in $\R^{n+1}$, 
the Gauss-Kronecker curvature is given by 
$$
K = \dfrac{(-1)^n \det\left(\frac{\pa^2 f}{\pa x_i\pa x_j}\right)}{\left[1 + \left(\frac{\pa f}{\pa x_1}\right)^2 
+ \cdots + \left(\frac{\pa f}{\pa x_n}\right)^2 \right]^{\frac{n+2}{n}}}. 
$$
Therefore, for a $C^2$-extension $\widetilde{S}$, $K$ must be continuous on $\widetilde{S}$. 
Thus, if $S$ is not flat in itself, 
then we have to impose just $C^1$-condition to the flat extensions $\widetilde{S}$. 

\

The efforts to solve the problem leads us to an insight on elementary differential 
geometry from singularity theory. 
We succeed the basic methods of geometric singularity theory(\cite{BG}\cite{Wall}). 
In fact we assume that the surface with boundary $(S, \gamma)$ is generic in this Introduction. 
However some of results hold for a surface {\it of finite type}: 
Regard it as a surface in the projective $3$-space 
$\R P^3$ and take its projective dual $S^{\vee}$ with boundary $\widehat{\gamma}$ 
in the dual space $\R P^{3*}$. Then the condition is that 
both $\gamma$ and $\hat{\gamma}$ are {\it of finite type} in the sense explained in 
\S \ref{Projective geometry of front-boundaries.}. 
Under the condition the tangent lines 
and the osculating planes to $\gamma$ (resp. $\widehat{\gamma}$) 
are well-defined. Generic surfaces with boundary are of finite type. 

A point $p \in \gamma$ is called an {\it osculating-tangent point} 
if the tangent plane $T_pS$ coincides with the osculating plane of $\gamma$, regarded as a space curve, at $p$. 
%

\bet
\label{solution of problem}
{\rm (The solution to generic flat extension problem).} 
Let $(S, \gamma)$ be a generic $C^\infty$ surface with boundary $\gamma$ in $\R^3$.  
Then $(S, \gamma)$ has a unique $C^1$ flat extension $\widetilde{S}$ 
locally across $\gamma$ 
near $p \in \gamma$ provided 
$p$ is not an osculating-tangent points for $(S, \gamma)$. 
\ent

\ber
{\rm
Let $g : S \to S^2$ be the Gauss mapping on $S$ in $\R^3$ (\cite{BGM}). 
Then the local uniqueness of flat extensions holds under the weaker condition that 
the spherical curve $g\vert_\gamma : \gamma \to S^2$ is immersive. 
}
\enr

In fact, to obtain the flat extension of $(S, \gamma)$ along the boundary $\gamma$, 
we take tangent planes to $S$ along $\gamma$ and 
take the {\it envelope} of the one-parameter family of tangent planes 
(See \S \ref{Projective geometry of front-boundaries.}. See also \cite{Thom}). 
We call it the {\it boundary-envelope} of $(S, \gamma)$. 
Then we have

\bet
\label{boundary-envelope}
For a generic $C^\infty$ surface $(S, \gamma)$ with boundary, 
the singularities of boundary-envelope of $(S, \gamma)$ 
are just cuspidal edges and swallowtails. 
\ent

\ber
{\rm 
The folded umbrella (or the cuspidal cross-cap) (\cite{Cleave}) does not appear 
as a generic singularity of boundary-envelope. It appears in a generic one parameter 
family of boundary-envelope (cf. Lemma \ref{key lemma} (2)). 
}
\enr

\bee
\label{example1}
{\rm
Let $S$ be a $C^\infty$ surface in $\R^3$ parametrised as 
$$
(x_1, x_2, x_3) = (t^2+u, \ t, \ t^3 + ut)
$$
with the parameters $t$ and $u$, the boundary $\gamma$ being given by 
$\{ u = 0\}$, namely, by 
$
\gamma(t) = (t^2, \ t, \ t^3). 
$
The osculating plane to $\gamma$ at $t = 0$ is given by $\{ x_3 = 0\}$ which is equal to 
the tangent plane of $S$ at the origin. Thus the origin is a osculating-tangent point 
of $(S, \gamma)$. 
Then the boundary-envelope of $(S, \gamma)$ is given by 
$$
(x, t) \mapsto (x_1, x_2, x_3) = (3t^2 - 2xt, \ x, \ -2t^3 + xt^2). 
$$
Its singular locus passes through the origin. 

}
\ene

\bee
\label{example}
{\rm 
Let $S$ be a $C^\infty$ surface in $\R^3$ parametrised as 
$$
(x_1, x_2, x_3) = (t+1, \ 4t^3 - 2t^2 - 2t + u, \ 3t^4 - t^3 - t^2 + ut)
$$
with the parameters $t$ and $u$, the boundary $\gamma$ being given by 
$\{ u = 0\}$, namely, by 
$$
\gamma(t) = (t+1, \ 4t^3 - 2t^2 - 2t, \ 3t^4 - t^3 - t^2). 
$$
Then the boundary-envelope of $(S, \gamma)$ is given by 
$$
(x, t) \mapsto (x_1, x_2, x_3) = (x, \ 4t^3 - 2xt, \ 3t^4 - xt^2). 
$$

\

In general a map-germ $(\R^2, 0) \to (\R^3, 0)$ is called a {\it swallowtail} (or of type $A_3$) if 
it is diffeomorphic, i.e. $C^\infty$ right-left equivalent, to the germ $(x, t) \mapsto (x, \ 4t^3 - 2xt, \ 3t^4 - xt^2)$ 
at $(0, 0)$. Moreover a map-germ $(\R^2, 0) \to (\R^3, 0)$ is called a {\it cuspidal edge} 
(or of type $A_2$) if 
it is $C^\infty$ right-left equivalent to the germ $(x, t) \mapsto (x, \ 3t^2 - 2xt, \ 2t^3 - xt^2)$ 
at $(0, 0)$. 

In our example, the cuspidal edge and the swallowtail singularities are realised by a flat surface, 
a $C^\infty$ surface which is flat outside the singular locus. 




Note that, in the above example, the dual surface $S^{\vee}$ is given by 
$$
(y_0, y_1, y_2) = (t^4 + u(t+1), \ t^2 + u, \ t)
$$
and its boundary $\widehat{\gamma}$ is given in $(y_0, y_1, y_2)$-space, by 
$$
\widehat{\gamma}(t) = (y_0(t), y_1(t), y_2(t)) = (t^4, \ t^2, \ t), 
$$
while $\widehat{\gamma}^*$ is given by 
$$
\widehat{\gamma}^*(t) = (x_1, x_2, x_3) = (6t^2, \ - 8t^3, \ - 3t^4), 
$$
for the notations which will be introduced in \S \ref{Projective geometry of front-boundaries.}. 
The singular locus of the boundary-envelope of $(S, \gamma)$ is given by 
$\widehat{\gamma}^*$. 
}
\ene

Motivated by this geometric method, we distinguish several \lq\lq landmarks", 
added to osculating-tangent points, 
on the boundary $\gamma$ for a generic surface $(S, \gamma)$: 
A {\it parabolic point} on the boundary $\gamma$ is a point on the intersection of the parabolic locus of $S$, 
the singular locus of the Gauss mapping $g : S \to S^2$ 
and $\gamma$ (\cite{BGM}). 
A point $p \in \gamma$ is called a {\it swallowtail-tangent point} 
if the tangent plane $T_pS$ contacts with the envelope at a swallowtail point of the envelope. 
It turns out that a point $t = t_1$ of the parametric boundary $\gamma$ is 
a swallowtail-tangent point if and only if, at $t = t_1$, 
the dual curve $(\widehat{\gamma})^*$ to the dual-boundary 
$\widehat{\gamma}$ is defined and the tangent line to the point $(\widehat{\gamma})^*$ at $t= t_1$ 
contains the swallowtail point of the envelope $(\widehat{\gamma})^{\vee}$. 

Parabolic points on $\gamma$ for $(S, \gamma)$ correspond to singular points of the dual $S^{\vee}$ 
on the dual-boundary $\widehat{\gamma}$. 

In Example \ref{example}, $\gamma$ has no osculating-tangent point nor 
parabolic point, but it has one swallowtail-tangent point at $(0, 0, 0)$. 

By Theorem \ref{solution of problem}, 
a generic surface with boundary $(S, \gamma)$ 
has a local flat extension 
across non-osculating-tangent points. In fact, at any 
osculating-tangent point, the singular locus of the boundary-envelope 
passes through the boundary at that point. 
See \S \ref{Projective geometry of front-boundaries.} for 
the exact classification of singularities of the local extension problem. 
Moreover a global obstruction occurs by singularities of the envelope, 
in particular, by self-intersection loci. 
Thus a swallowtail point of the envelope provides 
\lq\lq a global obstruction with local origin" 
for the flat extension problem.  
With this motivation, we characterise the osculating tangent points and 
the swallowtail tangent points in terms of Euclidean invariant of 
the surface-boundary $\gamma$ of $S$. 

To characterise these landmarks, 
we recall three fundamental invariants $\kappa_1, \kappa_2$ and $\kappa_3$ 
of the boundary $\gamma$ in \S \ref{Euclidean geometry of surface-boundaries.}. 
Actually 
$\kappa_1$ is the {\it geodesic curvature}, 
$\kappa_2$ is the {\it normal curvature} and 
$\kappa_3$ is the {\it geodesic torsion} of $\gamma$, up to sign. 
These three invariants are defined for any immersed space curve with 
a framing.

\ber
{\rm
The curvature $\kappa$ and the 
torsion $\tau$ of $\gamma$ as a space curve is related to 
$\kappa_1, \kappa_2$ and $\kappa_3$ by 
$$
\begin{array}{ccl}
\kappa & = & \sqrt{\kappa_1^2 + \kappa_2^2}, 
\vspace{0.2truecm}
\\
\tau & = & \kappa_3 + \left(\dfrac{\kappa}{\kappa_1}\right)
\left(\dfrac{\kappa_2}{\kappa}\right)' = \kappa_3 - 
\left(\dfrac{\kappa}{\kappa_2}\right)
\left(\dfrac{\kappa_1}{\kappa}\right)' = 
\kappa_3 + \dfrac{\kappa_1\kappa_2' - \kappa_2\kappa_1'}{\kappa_1^2 + \kappa_2^2}, 
\end{array}
$$
for the arc-length differential, provided $\kappa_1 \not= 0$ and 
$\kappa_2 \not= 0$. 
Note that the torsion $\tau$ of an immersed space curve is defined when the curvature $\kappa \not= 0$. 
Moreover it can be shown that, for any space curve $\gamma$ with curvature $\kappa$ and $\tau$, ($\kappa 
\not= 0$),  
and given any three functions $\kappa_1, \kappa_2$ and $\kappa_3$ 
on the curve satisfying the above relations, 
there exists a surface $S$ with boundary $\gamma$ such that the three invariants coincide with the 
given $\kappa_1, \kappa_2$ and $\kappa_3$. 
}
\enr

Then our generic characterisation is given by 

\bet
\label{osculating-tangent}
Let $(S, \gamma)$ be a generic $C^\infty$ surface with boundary 
in Euclidean three space $\R^3$. Then the osculating-tangent point on $\gamma$ 
is characterised by the condition $\kappa_2 = 0$. 
\ent 

Moreover, we show that there exists a characterisation of the swallowtail-tangent points 
in terms of $\kappa_1, \kappa_2, \kappa_3$ and their derivatives of order $\leq 3$. 
In fact we have 

\bet
\label{swallowtail-tangent}
{\rm (Euclidean generic characterisation of swallowtail-tangent)} 
Let $(S, \gamma)$ be a generic $C^\infty$ surface with boundary in Euclidean three space $\R^3$. 
A swallowtail-tangent point of $\gamma$ is characterised by the condition
\vspace{0.2truecm}
\\
{\rm (I)} $\kappa_2 \not= 0$, 
\vspace{0.2truecm}
\\
{\rm (II)} $\kappa_1^2\kappa_3(\kappa_2^2 + \kappa_3^2) 
+ \kappa_2(\kappa_2^2 + \kappa_3^2)\kappa_1' 
- 3\kappa_1\kappa_3^2\kappa_2' 
+ 3\kappa_1\kappa_2\kappa_3\kappa_3' 
+ 2\kappa_3(\kappa_2')^2 - 2\kappa_2\kappa_2'\kappa_3' 
- \kappa_2\kappa_3\kappa_2'' + \kappa_2^2\kappa_3'' \ = \ 0$, 
\vspace{0.2truecm}
\\
{\rm (III)} 
$
2\kappa_1\kappa_2^3(\kappa_1^2 + \kappa_2^2 + \kappa_3^2) 
+ 2\kappa_1\kappa_3(2\kappa_2^2 + \kappa_3^2)\kappa_1' 
+ (3\kappa_2^2 - 2\kappa_3^2)\kappa_1'\kappa_2' 
+ 5\kappa_2\kappa_3\kappa_1'\kappa_3' + 3\kappa_1\kappa_2(\kappa_3')^2 
+ \kappa_2(3\kappa_1\kappa_2 + \kappa_2^2 + \kappa_3^2)\kappa_1''
+ 3\{ \kappa_1(- \kappa_2^2 - \kappa_3^2 + \kappa_2\kappa_3)
+ 3(\kappa_3\kappa_2' - \kappa_2\kappa_3')\} \kappa_3'' 
+ \kappa_2(\kappa_2 - 2\kappa_3)\kappa_3'''
\not= \ 0$. 
\ent

\ber
{\rm 
\label{osculating-tangent-elliptic}
The existence of an osculating-tangent point on the boundary $\gamma$ 
depends on the geometry of the surface $S$ itself. 

For example, 
on an elliptic surface, there does not exist any osculating-tangent point. 
The surface is necessarily hyperbolic near an osculating-tangent point 
with $\kappa_2 = 0, \kappa_3 \not= 0$. 
}
\enr

We are interested in 
the interaction between singularity and geometry. 
In our topic of this paper, local geometry of surface-curve provides a global effect to the singularity of 
the envelope. In fact we give the exact formula for 
the distance between the swallowtail-tangent point 
on the surface-boundary
and the swallowtail point on the boundary-envelope (envelope-swallowtail) 
in terms of local geometric invariants of the boundary. (See also Proposition \ref{distance2}). 

\bep
\label{distance}
The distance $d$ between 
the swallowtail tangent point on the surface-boundary and 
the envelope-swallowtail is given by 
$$
d = \left\vert 
\dfrac{\kappa_2\sqrt{\kappa_2^2 + \kappa_3^2}}{\kappa_2(\kappa_3' + \kappa_1\kappa_2) 
+ \kappa_3(-\kappa_2' + \kappa_1\kappa_3)}
\right\vert. 
$$
\enp

\ber
{\rm 
If the denominator of the above formula 
vanishes, then the formula reads $d = \infty$, and, in fact, the envelope-swallowtail lies at infinity. 
If $\kappa_2 = 0$, then the formula reads $d = 0$, and, in fact, 
the non-generic coincidence of an osculating-tangent point and a swallowtail-tangent point occurs, 
and the envelope-swallowtail coincides with the swallowtail-tangent point. 
}
\enr


In \S \ref{Projective geometry of front-boundaries.}, 
we give the background for the basic results on projective duality of Legendre surfaces 
with boundaries (Theorems \ref{basic results} and \ref{envelope}). 
As a corollary we show Theorem \ref{solution of problem}. 
In \S \ref{Euclidean geometry of surface-boundaries.}, 
we show the Euclidean characterisations of osculating-tangent points 
and swallowtail-tangent points (Theorems \ref{osculating-tangent} and \ref{swallowtail-tangent}) 
and the distance formula 
(Proposition \ref{distance}) in more general setting: 
Our perspective though singularity theory extends the results on 
generic surface-boundaries 
to more general surface-boundaries. 
Lastly the local flat extension problem is solved naturally as a by-product of other results in this paper. 

A local geometry of surface-boundary causes a global effect to the singularities of 
boundary-envelope. Thus we provide examples of results on the interaction between singularity and geometry and between local and global. 

Apart from the flat extension problem, also there exist several extension problem: 
For instance we can consider the $C^1$ extension problem by a surface with $K = c$ for a non-zero constant $c$. 
Note that generically a surface of constant Gaussian curvature has only cuspidal edges and 
swallowtails as singularities (\cite{IM}). 

Some of the results in this paper have been announced in the monograph \cite{Ishikawa6}.

\section{Projective geometry on singularities of front-boundaries.}
\label{Projective geometry of front-boundaries.}

To study the existence, uniqueness and singularities of flat extensions by the geometric method, 
we recall several basic results on Legendre surfaces with boundaries in projective-contact framework 
(\cite{Bruce}). 

The projective duality between the projective $(n+1)$-space  $\R P^{n+1} = P(\R^{n+2})$ 
and the dual projective $(n+1)$-space $\R P^{{n+1}*} = P(\R^{{n+2}*})$ 
is given through the incidence manifold 
$$
I^{2n+1} = \{ ([X], [Y]) \in \R P^{n+1}\times \R P^{{n+1}*} 
\mid X \cdot Y  = 0 \}, 
$$
and projections $\pi_1 : I^{2n+1} \to \R P^{n+1}$ and $\pi_2 : I^{2n+1} \to \R P^{{n+1}*}$. 
The space $I$ is identified with the space $PT^*\R P^{n+1}$ of contact elements of $\R P^{n+1}$ 
and with $PT^*\R P^{{n+1}*}$ as well. See \cite{O.P.Scherbak1} for instance. 
It is endowed with the natural contact structure 
$$
D = \{ X \cdot dY = 0 \} = \{ dX \cdot Y = 0 \} \subset TI \cong T(PT^*\R P^{n+1}) \cong T(PT^*\R P^{{n+1}*}). 
$$

A $C^\infty$ hypersurface $S$ in $\R P^{n+1}$ lifts uniquely to the Legendre hypersurface $L$ in $I$ 
which is an integral submanifold to $D$: 
$$
L = \{ ([X], [Y]) \in I \mid [X] \in S, 
[Y] {\mbox{\rm \ determines\ }} T_{[X]}S {\mbox{\rm \ as a
projective hyperplane}} \}. 
$$
Then $L$ projects to $\R P^{{n+1}*}$ by $\pi_2$. 
The \lq\lq front" $S^{\vee} = \pi_2(L)$, 
as a parametrised hypersurface with singularities, is called the 
{\it projective dual} or {\it Legendre transform} of $S$ (\cite{AGV}). 

\

If we start with a surface $S$ with boundary $\gamma$ in $\R P^3$, $n=2$, 
then the Legendre lift $L$ also has the boundary $\Gamma$: 
$$
\Gamma = \{ ([X], [Y]) \in L \mid [X] \in \gamma \} = \pa L. 
$$
Then $L$ is a Legendre surface and $\Gamma$ is an integral curve in $I^5$ 
to the contact distribution $D$: 
$$
T\Gamma \subset TL \subset D \subset TI. 
$$

Now we have a Legendre surface with boundary $(L, \Gamma)$ 
in $I$ 
and two Legendre fibrations 
$\pi_1, \pi_2$: 
$$
\begin{array}{ccccccc}
(L, \Gamma) & \subset & PT^*\R P^3 & \cong & \quad I^5 & \cong & PT^*\R P^{3*} 
\vspace{0.2truecm}
\\
\downarrow &            &  
\downarrow  &
\hspace{-0.5truecm}
{\mbox{\footnotesize $\pi_1$}}\swarrow 
&   & \hspace{0.5truecm} \searrow {\mbox{\footnotesize $\pi_2$}} & 
\downarrow
\vspace{0.2truecm}
 \\
(S, \gamma) & \subset   & 
\R P^3 & & & & 
\R P^{3*} 
\end{array}
$$

We identify $\Gamma$ with the inclusion map $\Gamma \hookrightarrow I$. 
Then we get 
the triple of Legendre surfaces $(L, L_1, L_2)$ possibly with singularities in $I$:  
$$
L_1 = \{ ([X], [Y]) \mid [X] \in \pi_1(\Gamma), [Y] 
{\mbox{\rm \ is a tangent plane to }} \pi_1\circ \Gamma {\mbox{\rm \ at \ }} [X] \},  
$$ 
the projective conormal bundle of the space curve $\pi_1(\Gamma)$, and 
$$
L_2 = \{ ([X], [Y]) \mid [Y] \in \pi_2(\Gamma), [X] 
{\mbox{\rm \ is a tangent plane to }} \pi_2\circ \Gamma {\mbox{\rm \ at \ }} [Y] \} 
$$
the projective conormal bundle of the space curve $\pi_2(\Gamma)$. 

Moreover, the {\it dual surface} of the space curve $\pi_1(\Gamma)$ (resp. 
$\pi_2(\Gamma)$) is defined as the front $\pi_2(L_1)$ (resp. $\pi_1(L_2)$). 
Thus we have two fronts or frontal surfaces $\pi_1(L), \pi_1(L_2) \subset \R P^3$ and 
$\pi_2(L), \pi_2(L_1) \subset \R P^{3*}$ respectively. 

\

Starting from $C^\infty$ surface $(S, \gamma)$ with boundary in $\R P^3$, 
we have the Legendre-integral lifting  
$(L, \Gamma)$ in $I^5$. 
Then the {\it boundary-envelope} of $(S, \gamma)$ is defined by 
$\pi_1\vert_{L_2} : L_2 \to \R P^3$.  
Moreover $\pi_2\vert_{L_1}$ gives the boundary-envelope of the dual $(S^{\vee}, \widehat{\gamma})$. 

\ber
\label{What is envelope?}
{\rm 
In the above definition of \lq\lq projective conormal bundle" $L_2$, 
the interpretation of \lq\lq tangent plane" is not unique if $\pi_2\circ\Gamma$ is 
not an immersion. In this paper we mainly concern with the generic case where 
$\pi_2\circ\Gamma$ is an immersion (cf. Theorem \ref{basic results} (3). 
See also Remark \ref{What is envelope?2}). 
}
\enr

\

We call a pair of germs of fronts $(S, E)$, say $(\pi_1\vert_L, \pi_1\vert_{L_1})$, 
is of type $B_2$ (resp. $B_3$, $C_3$) if 
it is diffeomorphic, i.e. $C^\infty$ right-left equivalent, to the following local model as a multi-germ: 
\vspace{0.3truecm}
\\
$B_2 : \ 
(t_1, t_2) \mapsto (t_1, t_2, t_2^2), \pm t_2 \geq 0, \ (s_1, s_2) \mapsto (s_1, s_2, 0)$
\\
$C_3^{\pm} : \ 
(t_1, t_2) \mapsto (t_1, t_2, 0), \pm t_2 \geq 0, \ (s_1, s_2) \mapsto (s_1, -3s_2^2 - 2s_1s_2, 2s_2^3 + s_1s_2^2)$
\\
$B_3 : \ 
(t_1, t_2) \mapsto 
(t_1, -3t_2^2 - 2t_1t_2, 2t_2^3 + t_1t_2^2), \pm t_2 \geq 0, 
\ 
(s_1, s_2) \mapsto (s_1, s_2, 0)$.
\\
($\pm$ for $B_2, B_3$ give the same class.)

  \begin{center}
  \includegraphics[width=8truecm, height=7truecm, clip, bb=10 300 550 850]{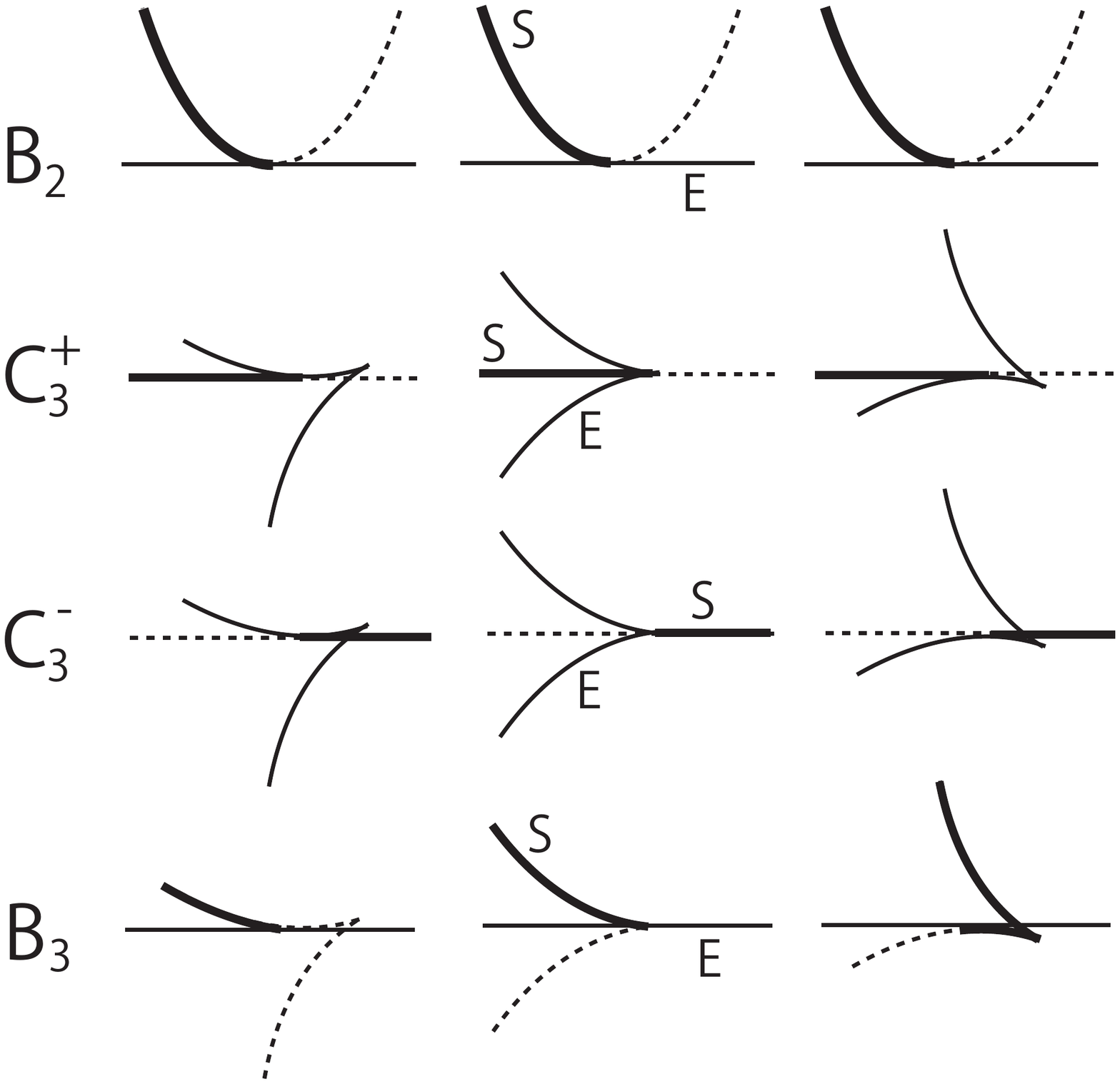}
\end{center}

Then we have the basic results: 

\bet
\label{basic results}
For a generic Legendre surface with boundary $(L, \Gamma)$ 
in the incident manifold $I^5 \cong PT^*\R P^3 \cong PT^*\R P^{3*}$ with respect to $C^\infty$ topology, 
we have 
\\
{\rm (1)} The singularities of $\pi_1\vert_L$ and $\pi_2\vert_L$ are 
just cuspidal edges and swallowtails. 
\\
{\rm (2)} The diffeomorphism types of the pair 
$(\pi_1\vert_L, \pi_1\vert_{L_2})$ {\rm (}resp. $(\pi_2\vert_L, \pi_2\vert_{L_1})${\rm )}
of germs at points on $\Gamma$ are given by $B_2, B_3$ and $C_3$. 
\\
{\rm (3)} Both $\pi_1\vert_{\Gamma}$ and $\pi_2\vert_{\Gamma}$ are generically 
immersed space curves in the sense of Scherbak {\rm (}\lq\lq Scherbak-generic"{\rm )} \cite{O.P.Scherbak1}, in $\R P^3$ and $\R P^{3*}$ respectively. 
Singularities of $\pi_1\vert_{L_2}$ and $\pi_2\vert_{L_1}$ are only cuspidal edges and swallowtails. 
\ent

\ber
{\rm 
We can show that the singular loci of $\pi_1\vert_L$ and $\pi_2\vert_L$, and $\Gamma$ 
are in general position in $L$ and 
moreover that the swallowtail points of $\pi_1\vert_L$ and $\pi_2\vert_L$ are 
not on the intersections of the above three curves. 
}
\enr

For the point (1), it is well-known that the stable front $\pi_1(L)$ 
has $A_\ell$-singularities $(\ell \leq 3)$ by Legendre singularity theory \cite{AGV}. 
The cuspidal edge singularity is called of type 
$A_2$ and the swallowtail singularity is called of type $A_3$, while $A_1$ means regular.  
For the point (2), 
it is well-known that 
the stable front with boundary $(\pi_1(L), \pi_1(L_2))$ has $B_\ell$ or $C_\ell$-singularity 
$(\ell \leq 3)$ by the theory of boundary singularities; 
we know the diffeomorphism types of stable fronts with boundary \cite{Arnold1}\cite{AGV}. 
See also \cite{I.G.Scherbak1}\cite{I.G.Scherbak2}\cite{DD}\cite{Tsukada}. 
Moreover, for the point (2), we remark that, the duality of boundary singularities found by 
I.G. Scherbak, the \lq\lq Scherbak duality" (\cite{I.G.Scherbak1}\cite{I.G.Scherbak2}) 
are realised via Legendre duality in our geometrical situation: 
The $C_3$-singularity appears at an osculating-tangent point on $\gamma$ in $\R P^3$ 
and $B_3$-singularity appears at a point in $\R P^{3*}$ 
corresponding to a parabolic point on $\gamma$. 

These basic results are proved by the standard methods in singularity theory: Here we use  
Legendre-integral version of relative transversality theorem \cite{Ishikawa1}\cite{Ishikawa5} 
to make assure ourselves. 
A Legendre immersion $i : (L, \Gamma) \to I^5$ is approximated by $i' : 
(L, \Gamma) \to I^5$ such that the $r$-jet extension $j^ri'$ is transverse 
to given a finite family of submanifolds in the isotropic jet space $J_{\rm{int}}^r(L, I)$ and 
$(j^ri')\vert_\Gamma : \Gamma \to J^r_{\rm{int}}(L, \Gamma; I, I)$ is 
transverse to given a finite family of submanifolds in the relative isotropic jet space 
$J^r_{\rm{int}}(L, \Gamma; I, I)$ which is a fibration over $\Gamma\times I$ (\cite{Ishikawa1}). 
Moreover $j^r(i'\vert_\Gamma)$ is 
transverse to given finite family of submanifolds in $J_{\rm{int}}^r(\Gamma, I)$.  
We will give a proof of 
the Legendre (or integral) transversality theorem, because it seems to be never explicitly given. 
%

\bet
\label{Integral transversality theorem}
{\rm (Integral transversality theorem \cite{Ishikawa3}\cite{Ishikawa5})} 
Let $(I^{2n+1}, D)$ be a $(2n+1)$-dimensional contact manifold, $M^m$ an $m$-dimensional 
manifold ($m \leq n)$ and $f : M \to I$ an integral immersion to the contact structure $D \subset 
TI$. 
Let $r \in \NN$ and $Q_\lambda (\lambda \in \Lambda)$ a finite family of submanifolds 
of $J^r_{\rm{int}}(M, I)$. 
Then $f$ is approximated, in the Whitney $C^\infty$ topology, by an integral immersion 
$f' : M \to I$ such that the $r$-jet extension $j^rf' : M \to J^r_{\rm{int}}(M, I)$ is transverse to 
all $Q_\lambda (\lambda \in \Lambda)$. 
\ent

\Proof 
First recall that the space of integral immersion-jets 
$J^r_{\rm{int}}(M, I)$ is a submanifold of $J^r(M, I)$ (\cite{Ishikawa3}). 
Then we follow the standard construction of \cite{Mather} in the integral context: 
Suppose, near each point $p \in M$ and $f(p) \in I$, $f$ is represented as 
$$
(t_1, \dots, t_m) \mapsto (t_1, \dots, t_m, 0, \dots, 0)
$$
by a local coordinates of $(M, p)$ and a local Darboux coordinates of $(I, f(p))$. 
Denote by $P(m, \ell; k)$ the space of polynomial mappings $\R^m 
\to \R^\ell$ of degree $\leq k$. 
Let $E$ be a neighbourhood of $(0, \id_{\R^m})$ of $P(m, 1; r+1)\times P(m, m; r)$. 
Choose a $C^\infty$ function $\rho : 
\R^m \to [0, 1]$ with a compact support. 
For $(S, \sigma)$, set 
$$
\varphi(S, \sigma)(t) = \left(t, \rho(t)\dfrac{\pa S}{\pa t}\right)\circ \sigma(t), 
$$
and extend it to an integral immersion $\varphi(S, \sigma) : M \to I$. 
Then $\Phi : E\times I \to J^r_{\rm{int}}(M, I)$ 
defined by $\Phi(S, \sigma, t) = j^r(\varphi(S, \sigma))(t))$ 
is a submersion at $(0, \id, p)$ and transverse to $Q_\lambda$ locally. 
Then the result follows by Sard's theorem. 
\QED

\ber
\label{RIT}
{\rm 
The relative version of Theorem \ref{Integral transversality theorem} 
is also valid, similarly to the construction of \cite{Ishikawa1}: 
Let $N^\ell \subset M^m$ be a submanifold and $r \in \NN$. 
Then we consider the relative integral jet space $J^r_{\rm{int}}(M, N, I, I)$ 
fibered over $N\times I$ with the fibre $J^r_{\rm{int}}(m, 2n+1)$,  
the space of jets of integral immersion-germs 
$(\R^m, 0) \to (\R^{2n+1}, 0)$ to a local model $\R^{2n+1}$ of the contact space. 
Let 
$Q_\lambda (\lambda \in \Lambda)$ be a finite family of submanifolds 
of $J^r_{\rm{int}}(M, I)$, $R_{\lambda'} (\lambda' \in \Lambda')$ a countable family of submanifolds 
of $J^r_{\rm{int}}(M, N, I, I)$ and 
$P_{\lambda''} (\lambda'' \in \Lambda'')$ 
a countable family of submanifolds of $J^r_{\rm{int}}(N, I)$. 
Then any integral immersion $f : M \to I$ Let $r \in \NN$ is approximated, in the Whitney $C^\infty$ topology, by an integral immersion 
$f' : M \to I$ such that the $r$-jet extension $j^rf' : M \to J^r_{\rm{int}}(M, I)$ is transverse to 
all $Q_\lambda (\lambda \in \Lambda)$, 
$j^rf'\vert_N : N \to J^r_{\rm{int}}(M, N, I, I)$ is transverse to 
all $R_\lambda' (\lambda' \in \Lambda')$ and 
$j^r(f'\vert_N) : N \to J^r_{\rm{int}}(N, I)$ is transverse to 
all $P_\lambda'' (\lambda'' \in \Lambda'')$. 
}
\enr

The genericity for the points (1) (2) is described in terms of {\it generating families}: 
In the affine open subset $U\times V = \{ X_0 \not= 0, Y_3 \not= 0\}$ of $\R P^3 \times \R P^{3*}$, 
we set $x_i = -X_i/X_0, y_j = -Y_{3-j}/Y_3, (1 \leq i, j \leq 3)$, and 
$$
F(x_1, x_2, x_3, y_1, y_2, y_3) = - y_3 + x_1y_2 + x_2y_1 - x_3. 
$$
Then $I \cap (U\times V)$ is defined by $F = 0$. 

Let $L = \{ (x_1(u, v), x_2(u, v), x_3(u, v), y_1(u, v), y_2(u, v), y_3(u, v)) \}$ 
be a Legendre surface in $I \cap (U\times V)$ parametrised 
by $(u, v) \in \R^2$. Then we have two families of functions 
$F_2, F_1 : \R^2\times \R^3 \to \R$, 
$$
\begin{array}{ccc}
F_2(u, v; y_1, y_2, y_3) & = & F(x_1(u, v), x_2(u, v), x_3(u, v), y_1, y_2, y_3), 
\vspace{0.2truecm}
\\
F_1(u, v; x_1, x_2, x_3) & = & F(x_1, x_2, x_3, y_1(u, v), y_2(u, v), y_3(u, v)). 
\end{array}
$$
Then $F_1$ (resp. $F_2$) is a generating family for $\pi_1\vert L$ (resp. $\pi_2\vert L$). 
Similarly, for an integral curve $\Gamma = 
\{ (x_1(t), x_2(t), x_3(t), y_1(t), y_2(t), y_3(t)) \}$, we set 
$$
\begin{array}{ccc}
G_2(t; y_1, y_2, y_3) & = & F(x_1(t), x_2(t), x_3(t), y_1, y_2, y_3), 
\vspace{0.2truecm}
\\
G_1(t; x_1, x_2, x_3) & = & F(x_1, x_2, x_3, y_1(t), y_2(t), y_3(t)). 
\end{array}
$$
Then $G_1$ (resp. $G_2$) is a generating family for $\pi_1\vert_{L_2}$ (resp. $\pi_2\vert_{L_1}$). 
Note that 
$
G_i = F_i\vert_{\gamma\times\R^3}, i = 1, 2. 
$
The local singularities of pair of fronts $(\pi_1\vert_L, \pi_1\vert_{L_2})$ 
(resp. $(\pi_2\vert_L, \pi_2\vert_{L_1})$) are represented, via the analysis of generating families, 
as strata in the integral jet spaces. 
Let us make clear the relation of transversality in Legendre jet space and that 
in the jet space of generation functions: 
We use the following basic method to show genericity: 

\bep
\label{transversality1}
Let $I^{2n+1}$ be a $(2n+1)$-dimensional contact manifold, 
$\pi : I \to B^{n+1}$ a Legendre fibration, and $N$ an $n$-manifold. 
Let $f : N \to I$ be a Legendre immersion, 
$u_0 \in N$ and $F : (N\times \Lambda, (u_0, \lambda_0)) \to \R$  a 
generating family of $f : (N, u_0) \to I$. 
Then we have 
\\
{\rm (1)} $f : (N, u_0) \to I$ is Legendre stable if and only if 
$F$ is ${\mathcal K}$-stable unfolding of 
$F\vert_{N\times\{\lambda_0\}}$. 
\\
{\rm (2)} $f : (N, u_0) \to I$ is Legendre stable if and only if $j^{n+1}f : 
(N, u_0) \to J^{n+1}_{\rm{int}}(N, I)$ is transversal to the Legendre orbit of $j^{n+1}f(u_0)$. 
\\
{\rm (3)} $F$ is ${\mathcal K}$-stable unfolding of $F\vert_{N\times\{\lambda_0\}}$ 
if and only if $j^{n+2}_1F : (\Lambda, \lambda_0) \to J^{n+2}(N, \R)$, defined by 
$j^{n+2}_1F(\lambda) = j^{n+2}(F\vert_{N\times\{\lambda\}})(u_0)$ is transverse to 
${\mathcal K}$-orbit of $j^{n+2}(F\vert_{N\times\{\lambda_0\}})(u_0)$. 
\enp

\ber
{\rm 
It is known that any germ of Legendre stable Legendre immersion $f : (N^n, u_0) 
\to I^{2n+1}$ is $(n+1)$-determined among Legendre immersion-germs(\cite{Ishikawa5}). 
Moreover, for its generating family $F : (N\times\Lambda, (u_0, \lambda_0)) \to \R$,  
$F\vert_{N\times\{\lambda_0\}} : (N\times\{\lambda_0\}, (u_0, \lambda_0)) \to \R$ 
is $(n+2)$-determined. Note that the $(k+1)$-jet of $F\vert_{N\times\{\lambda_0\}}$ is determined by the $k$-jet of $f$ from the integrality condition. 
}
\enr

Also we use the relative version: 

\bep
\label{transversality2}
Let $I^{2n+1}$ be a $(2n+1)$-dimensional contact manifold, 
$\pi : I \to B^{n+1}$ a Legendre fibration, and 
$(N, \pa N)$ an $n$-manifold with boundary. 
Let $f : (N, \pa N) \to I$ be a Legendre immersion, 
$u_0 \in \pa N$ and $F : (N\times \Lambda, (u_0, \lambda_0)) \to \R$  
a generating family of $f : (N, u_0) \to M$. 
Then we have 
\\
{\rm (1)} $f : (N, \pa N, u_0) \to I$ is Legendre stable if and only if 
$F$ is ${\mathcal K_b}$-stable unfolding of 
$F\vert_{N\times\{\lambda_0\}}$. 
\\
{\rm (2)} $f : (N, \pa N, u_0) \to I$ is Legendre stable if and only if $j^{n+1}f\vert_{\pa N} : 
(\pa N, u_0) \to J^{n+1}_{\rm{int}}(N, \pa N; I, I)$ 
is transversal to Legendre orbit of $j^{n+1}f(u_0)$. 
\\
{\rm (3)} $F$ is ${\mathcal K_b}$-stable unfolding of $F\vert_{N\times\{\lambda_0\}}$ 
if and only if $j^{n+2}_1F : (\Lambda, \lambda_0) \to J^{n+2}(N, \pa N; \R, \R)$, defined by 
$j^{n+2}_1F(\lambda) = j^{n+2}(F\vert_{N\times\{\lambda\}})(u_0)$ is transverse to 
${\mathcal K}$-orbit of $j^{n+2}(F\vert_{N\times\{\lambda_0 \}})(u_0)$. 
\enp

In the above Proposition, ${\mathcal K_b}$-equivalence means boundary ${\mathcal K}$-equivalence 
The points (1)(3) are basic results in (boundary) singularity theory(\cite{AGV}). 
The point (2) follows the infinitesimal characterisation of Lagrange stability. 
See \cite{Ishikawa3}.

\ 

For the point (3) of Theorem \ref{basic results}, 
we have to know more information on the projective geometry of boundaries, 
$\gamma = \pi_1(\Gamma)$ and $\widehat{\gamma} = \pi_2(\Gamma)$. 
We write $\gamma = \pi_1(\Gamma)$ and $\widehat{\gamma} = \pi_2(\Gamma)$, 
and call $\widehat{\gamma}$ the {\it dual-boundary} to $\gamma$. 

To show the point (3), we recall some projective geometry-singularity in three space: 
We use, to a space curve $c$ in $\R P^3$ (resp. in $\R P^{3*}$), 
the notions of the dual curve $c^*$ and the dual surface 
$c^{\vee}$ in $\R P^{3*}$ (resp. $\R P^3$). 
Note that the dual-boundary $\widehat{c}$ is different from the dual curve $c^*$ 
to $c$ and it is defined only when $c$ is regarded as a surface-curve or a {\it 
framed} curve. 

A $C^\infty$ space curve $\gamma : \R \to \R P^3$ is called 
{\it of finite type} at $t = t_0 \in \R$, 
if for each system of affine coordinates in $\R P^3$ near $\gamma(t_0)$, 
the $3\times\infty$ matrix 
$$
(\gamma'(t_0), \gamma''(t_0), \dots, \gamma^{(r)}(t_0), \dots)
$$
is of rank $3$. Introduce the $3\times r$-matrix 
$$
A_r(t) = (\gamma'(t), \gamma''(t), \dots, \gamma^{(r)}(t)). 
$$
Then the type $(a_1, a_2, a_3)$ of $\gamma$ at $t = t_0$ is define by 
$$
\begin{array}{c}
a_1 = \min \{ r \mid \rank A_r(t_0) = 1\}, \quad a_2 = \min \{ r \mid \rank A_r(t_0) = 2\}, 
\vspace{0.2truecm} 
\\
a_3 = \min \{ r \mid \rank A_r(t_0) = 3\}. 
\end{array}
$$
Remark that $a_1, a_2, a_3$ are positive integers 
with $a_1 < a_2 < a_3$ and that, for 
some system of affine coordinates centred at $\gamma(t_0)$, 
$\gamma$ is expressed as 
$$
\left\{
\begin{array}{rcl}
X_1(t) & = & (t - t_0)^{a_1} + o((t - t_0)^{a_1}), 
\vspace{0.2truecm}\\
X_2(t) & = & (t - t_0)^{a_2} + o((t - t_0)^{a_2}), \vspace{0.2truecm}\\
X_3(t) & = & (t - t_0)^{a_3} + o((t - t_0)^{a_3}). 
\end{array}
\right.
$$
A point of $\gamma$ of type $(1, 2, 3)$ is called an ordinary point. 
Otherwise, it is called a special point of $\gamma$.  
Special points 
are isolated on $\R$ for a space curve of finite type. 

\bel
{\rm (O.P. Scherbak \cite{O.P.Scherbak1}):}  
A generic space curve $\gamma$ in $\R P^3$ is of type $(1, 2, 3)$ 
or $(1, 2, 4)$ at each point. 
\enl

\Proof
Consider the $3$-jet space 
$$
J^3(\R, \R P^3) = \{ j^3\gamma(t_0) \mid 
\gamma = (X_0(t), X_1(t), X_2(t), X_3(t))  : (\R, t_0) \to \R P^3 \} 
$$ 
of curves in $\R P^3$. 
Set 
$$
\Sigma = \{ j^3\gamma(t_0) \mid \det(\gamma(t_0), \gamma'(t_0), \gamma''(t_0), \gamma'''(t_0)) = 0 \}. 
$$
The conditions are independent of the choice of homogeneous coordinates of $\gamma$. 

Then $\Sigma$ is a fibration over $\R\times\R P^3$ 
whose fibre is an algebraic hypersurface 
in the jet space $J^3(1, 3)$. 
A map-germ $\gamma : (\R, t_0) \to \R P^3$ is of type $(1, 2, 3)$ (resp, $(1, 2, 4)$) 
if and only if $ j^3\gamma(t_0) \not\in \Sigma$  
(resp. $j^3\gamma(t_0) \in \Sigma$ and 
$j^3\gamma : (\R, t_0) \to J^3(\R, \R P^3)$
is transverse to $\Sigma$). 
Therefore, by the transversality theorem, we have the result. 
\QED

We call a curve {\it Scherbak-generic} if it is of finite type of type $(1, 2, 3)$ or $(1, 2, 4)$ 
at any point. 

\

The osculating planes to a space curve $\gamma$ of finite type 
form a dual curve $\gamma^*$ of the curve $\gamma$ in the dual space. 

\bel
{\rm (Duality Theorem, Arnol'd, Scherbak \cite{O.P.Scherbak1}): }
\\
{\rm (1)} The dual curve $\gamma^*$ to a curve-germ $\gamma$ 
of finite type $(a_1, a_2, a_3)$ is 
a curve-germ of finite type $(a_3 - a_2, a_3 - a_1, a_3)$. 
\\
{\rm (2)} 
The dual surface to a curve-germ $\gamma$ of finite type is 
the tangent developable
of the dual curve $\gamma^*$ of $\gamma$.  
\enl

The {\it tangent developable} of $\gamma$ is a surface ruled by tangent lines to $\gamma$ 
(\cite{Cleave}\cite{Mond1}\cite{Mond2}\cite{O.P.Scherbak1}\cite{O.P.Scherbak2}\cite{Ishikawa4}). 

\ber
\label{What is envelope?2} 
{\rm 
The notion of dual surface depends on the notion of tangency (Remark \ref{What is envelope?}). 
For the curves of finite type, the notion of tangent line is well defined. 
Therefore if $\pi_1\circ\Gamma$ and $\pi_2\circ\Gamma$ 
are both of finite type, then 
both $L_1, L_2$ are well-defined, so are both $\pi_1\vert_{L_2}$ and $\pi_2\vert_{L_1}$. 
Thus the notion of boundary-envelope is well-defined. 
Also note that if we start from the generating family to define the boundary-envelope, 
we get the \lq\lq extended" envelope: To each singular point of $\pi_2\circ\Gamma$ 
the hyperplane in $\R P^3$ which corresponds to it is added to the original envelope $\pi_1(L_2)$. 
}
\enr

\bel
\label{dual-cp-sw}
If $\gamma$ is of type $(1, 2, 3)$, then $\gamma^*$ is of type $(1, 2, 3)$, and 
the dual surface is diffeomorphic to the cuspidal edge. 
If $\gamma$ is of type $(1, 2, 4)$, then $\gamma^*$ is of type $(2, 3, 4)$, and 
the dual surface is diffeomorphic to the swallowtail. 
\enl

For the proof, consult the survey paper \cite{Ishikawa4} 
on the singularities of tangent developables. 
We also remark 

\bel
\label{char}
The dual surface of a space curve-germ $\gamma$ of finite type is diffeomorphic to the cuspidal edge 
(resp. the swallowtail) 
if and only if the type of $\gamma$ is equal to $(1, 2, 3)$ {\rm (}resp. $(1, 2, 4)${\rm )}. 
\enl

Note that the type of $\widehat{\gamma}^*$ is $(1, 2, 3)$ (resp. 
$(2, 3, 4)$) if and only if $\widehat{\gamma}$ is of type $(1, 2, 3)$ (resp. 
$(1, 2, 4)$). 
Then Lemma \ref{char} follows from 
the following general result which does not stated in \cite{Ishikawa4}: 

\bep
Let $\gamma, \gamma'$ be space curve-germs of finite types. 
If their tangent developables 
are 
diffeomorphic, then their types 
coincide. 
\enp

\Proof
Let ${\mbox{\rm{type}}}(\gamma) = (m, m+s, m+s+r)$. Then diffeomorphism-class of the 
tangent developable of $\gamma$ is 
given by ${\mbox{\rm{dev}}}(\gamma): (\R^2, 0) \to (\R^3, 0)$ 
$$
x_1 = x, \ x_2 = t^{s+m} + \cdots + x(t^s + \cdots), \ 
x_3 = t^{r+s+m} + \cdots + x(ct^{r+s} + \cdots), 
$$
where $(x, t)$ is a system of parameters, $\cdots$ means higher order terms in $t$, 
and $c$ is a non-zero constant (\cite{Ishikawa4}\cite{Ishikawa2}). 
Suppose ${\mbox{\rm{dev}}}(\gamma)$ and ${\mbox{\rm{dev}}}(\gamma')$ are diffeomorphic by 
diffeomorphism-germs $\sigma : (\R^2, 0) \to (\R^2, 0)$ and $\tau: (\R^3, 0) \to (\R^3, 0)$, 
and the type of $\gamma'$ is $(m', m'+s', m'+s'+r')$. 
In general ${\mbox{\rm{dev}}}(\gamma)$ has singularity always along the original space curve $\gamma$, 
$\{ x = 0\}$, 
and along the tangent line to $\gamma$ at the origin $\{ t = 0\}$ when $s \geq 2$. 
Furthermore ${\mbox{\rm{dev}}}(\gamma)$ has the cuspidal edge 
singularity along $x = 0, t \not= 0$, while 
it has singularity along $\{ t = 0, x \not= 0\}$ if and only if $s = 2, r = 1$. 
On the other hand the curve $\gamma$ itself is singular if and only if $m \geq 2$. 
Therefore if the type is not equal to $(1, 3, 4)$, then the diffeomorphism $\sigma$ preserves 
$\{ x = 0\}$. Then $\sigma$ and $\tau$ have some restrictions: 
The first component of $\sigma$ is of form $x\rho(x, t)$, $\rho(0, 0) \not= 0$. 
The linear term of $\tau$ preserves the plane $\{ x_1 = 0\}$. 
Therefore, by the order comparison on $t$, we see that 
$s + m = s' + m',  r+s+m = r'+s'+m'$. 
Moreover, restricting the equivalence on $\gamma$ (and $\gamma'$), we see 
$m = m'$. Hence we have $(m, m+s, m+s+r) = (m', m'+s', m'+s'+r')$. 
\QED

\

A $C^\infty$ surface $(S, \gamma)$ with boundary is called {\it of finite type} if the boundary $\gamma$ 
and the dual-boundary $\widehat{\gamma}$ are both of finite type. 
Note that generic surfaces are of finite type (Lemma \ref{generic is finite}). 

From the above argument, in particular we have 

\bel
\label{envelope}
If $(S, \gamma)$ is of finite type, then the boundary-envelope of $(S, \gamma)$ 
is the dual surface $(\widehat{\gamma})^{\vee}$ of the dual-boundary $\widehat{\gamma}$. 
The boundary-envelope is the tangent developable to the dual curve $(\gamma^\vee)^*$ to 
the dual-boundary $\gamma^\vee$. 
Moreover, if $(S, \gamma)$ 
is generic, then 
there are only cuspidal edge singularities and swallowtail singularities on the boundary-envelope 
$\pi_1\vert_{L_2}$. 
\enl

\ber
To investigate the global flat extension problem, we need the global study on singularities 
of tangent developables. For this subject, see \cite{NS}. 
\enr

The following lemma is also a key for the theory:

\bel
\label{key lemma}
Let $I$ be a $(2n+1)$-dimensional contact manifold, 
$\pi : I \to B$ a Legendre fibration over an $(n + 1)$-dimensional manifold $B$, and 
$k\geq 1$. 
\\
{\rm (1)} 
Define $\Pi : J_{\rm{int}}^k(\R, I) \to J^k(\R, B)$ 
by $\Pi(j^k\Gamma(t_0)) = j^k(\pi\circ\Gamma)(t_0)$ 
for any integral curve-germ $\Gamma : (\R, t_0) \to I$. 
Then $\Pi$ is a submersion at $j^k\Gamma(t_0)$ if $\pi\circ\Gamma$ is 
an immersion at $t_0$
\\
{\rm (2)} 
The set $\Sigma = \{ j^k\Gamma(t_0) \in J_{\rm{int}}^k(\R, I) \mid 
\pi\circ\Gamma {\mbox{\rm  \ is not an immersion at\ }} t_0\}$ is of codimension 
$n$ in $J_{\rm{int}}^k(\R, I)$. 
\enl

\ber
\label{generic is finite}
{\rm 
By Lemma \ref{key lemma} ($n = 2$), we have the following: 
Let $\Pi_1 : J_{\rm{int}}^k(\R, I^5) \to J^k(\R, \R P^3)$ 
(resp. $\Pi_2 : J_{\rm{int}}^k(\R, I^5) \to J^k(\R, \R P^{3*})$) be the mapping 
induced by the Legendre fibration $\pi_1 : I \to \R P^3$ (resp. $\pi_2 : I \to \R P^{3*}$). 
Then the set $\Sigma_1$ (resp. $\Sigma_2$) of jets with singularity after the 
projection $\pi_1$ (resp. $\pi_2$) is of codimension $2$ in $J_{\rm{int}}^k(\R, I^5)$. 
Moreover $\Pi_1 : J_{\rm{int}}^k(\R, I^5) \setminus \Sigma_1 \to 
J^k(\R, \R P^3) \setminus \Pi_1(\Sigma_1)$ (resp. $\Pi_2 
: J_{\rm{int}}^k(\R, I^5) \setminus \Sigma_2 \to 
J^k(\R, \R P^{3*}) \setminus \Pi_2(\Sigma_2)$) is a submersion. 
}
\enr

\noindent
{\it Proof of Lemma \ref{key lemma}.} 
Take Darboux coordinates 
$x_1, \dots, x_n, z, p_1, \dots, p_n$ of $I$ around $\Gamma(t_0)$ and 
$x_1, \dots, x_n, z$ of $B$ around $\pi\circ\Gamma(t_0)$ so that 
the contact structure is given by $dz - (p_1 dx_1 + \cdots + p_n dx_n) = 0$ and 
$\pi$ is given by $(x_1, \dots, x_n, z, p_1, \dots, p_n) \mapsto (x_1, \dots, x_n, z)$. 
\\
(1) 
Set 
$
\Gamma(t) = (x_1(t), \dots, x_n(t), z(t), p_1(t), \dots, p_n(t))
$ and 
suppose $\pi\circ\Gamma$ is an immersion at $t_0$. 
Without loss of generality, we suppose $\dot{x}_1(t_0) \not= 0$. 
Take any deformation $c(t, s) = (X_1(t, s), \dots, X_n(t, s), Z(t, s))$ of $\pi\circ\Gamma(t)$ at $s = 0$. 
Note that $\dot{z} = p_1\dot{x}_1 + \cdots + p_n\dot{x}_n$. 
Therefore $p_1(t) = \dfrac{\dot{z}(t)}{\dot{x}_1(t)} - p_2\dfrac{\dot{x}_2(t)}{\dot{x}_1(t)} - 
\dots - p_n\dfrac{\dot{x}_n(t)}{\dot{x}_1(t)}$, near $t = t_0$. 
We set 
$$
\begin{array}{ccl}
P_1(t, s) & := & \dfrac{\dot{Z}(t, s)}{\dot{X}_1(t, s)} - p_2(t)\dfrac{\dot{X}_2(t, s)}{\dot{X}_1(t, s)} - 
\dots - p_n(t)\dfrac{\dot{X}_n(t, s)}{\dot{X}_1(t, s)}, 
\vspace{0.2truecm}
\\
P_i(t, s) & := & p_i(t), (i = 2, \dots, n), 
\end{array}
$$
near $(t, s) = (t_0, 0)$. Here $\dot{Z}(t, s)$ means the derivative by $t$. 
Then we get the integral deformation 
$$
C(t, s) = (X_1(t, s), \dots, X_n(t, s), Z(t, s), P_1(t, s), \dots, P_n(t, s))
$$ 
of $\Gamma(t)$ at $s = 0$, which satisfies 
$\pi(C(t, s)) = c(t, s)$. This show that any curve starting at $j^k(\pi\circ \Gamma)(t_0)$ 
in $J^k(\R, B)$ lifts to a curve starting at $j^k\Gamma(t_0)$ 
in $J_{\rm{int}}^k(\R, I)$. Therefore $\Pi$ is a submersion at $j^k\Gamma(t_0)$. 
\\
To see (2), first remark that $J_{\rm{int}}^k(\R, I)$ has local coordinates 
$$
x_1^{(i)}, \dots, x_n^{(i)}, z, p_1^{(i)}, \dots, p_n^{(i)}, \quad (0 \leq i \leq k), 
$$
because $z^{(i)} (2 \leq i \leq k)$ are written by these coordinates from the integrality condition
$\dot{z} = p_1\dot{x}_1 + \cdots + p_n\dot{x}_n$. 
Then $\Sigma$ is defined exactly by $x_1' = \cdots = x_n' = 0$. 
\QED


\

\noindent
{\it 
Proof of Theorem \ref{basic results}: 
}
As is mentioned above, 
we prove Theorem \ref{basic results} 
using relative version of Theorem \ref{Integral transversality theorem}
instead of the ordinary transversality theorem. 
In fact, we consider three kinds of transverslities: 
Transversality in $J^k_{\rm{int}}(\R^2, I^5)$, that in 
$J^k_{\rm{int}}(\R^2, \R; I^5, I^5)$ and that in $J^k_{\rm{int}}(\R, I^5)$. 
Note that the relative jet space $J^k_{\rm{int}}(\R^2, \R; I^5, I^5)$ is fibered over 
$\R\times I$ with fiber $J^k_{\rm{int}}(2, 5)$, the space of jets of integral immersions 
$(\R^2, 0) \to (\R^5, 0)$ to a local model $\R^5$ of the contact space, which is 
the fibre also for $J^k_{\rm{int}}(\R^2, I)$. 
However we consider the group action on $J^k_{\rm{int}}(2, 5)$ 
for $J^k_{\rm{int}}(\R^2, I)$ (resp. $J^k_{\rm{int}}(\R^2, \R; I, I)$) by 
diffeomorphisms on $(\R^2, 0)$ (res. 
by relative diffeomorphisms on $(\R^2, \R)$ ) 
and fiber-preserving contactomorphisms 
on $(\R^5, 0)$ with a local model $(\R^5, 0) \to (\R^3)$ of Legendre fibration. 
We take $k$ sufficiently large. Actually it is enough to take $k \geq 3$ in our case. 
We use Propositions \ref{transversality1} and \ref{transversality2}. 
In $J^k_{\rm{int}}(\R^2, I)$, we see that the complement to the union of  
$A_\ell$-orbits $(\ell \leq 3)$ is of codimension $3$ in the jet space 
of Legendre immersions $J^k_{\rm{int}}(\R^2, I)$. 
Then we have (1) by Theorem \ref{Integral transversality theorem}. 
Moreover, in $J^k_{\rm{int}}(\R^2, \R; I, I)$, 
the complement to the union of $B_\ell$ and $C_\ell$-orbits $(\ell \leq 3)$ 
is of codimension $2$ in $J^k_{\rm{int}}(\R^2, \R; I, I)$ along boundary 
(cf. \cite{Arnold1} Theorem 1, Remark 1). 
Thus, by the relative integral transversality theorem (Remark \ref{RIT}), 
we have (2). 
In $J^k_{\rm{int}}(\R, I)$, by Lemma \ref{key lemma} and Remark \ref{generic is finite}, 
we see that the complement to 
the jets of integral curves $\Gamma$ such that $\pi_1\circ \Gamma$ (resp. 
$\pi_2\circ \Gamma$) is Scherbak-generic, is of codimension $2$. 
Therefore, by Theorem \ref{Integral transversality theorem}, we have that, 
for a generic integral curve $\Gamma$ in $I$, both $\pi_1\vert_\Gamma$ and 
$\pi_2\vert_\Gamma$ are Scherbak-generic. 
Therefore by Lemma \ref{dual-cp-sw}, we have (3). 
\QED


\section{Euclidean geometry of surface-boundaries.}
\label{Euclidean geometry of surface-boundaries.}

The fundamental construction to observe such characterisations 
as Theorems \ref{osculating-tangent} and \ref{swallowtail-tangent} is as follows: 

The unit tangent bundle 
$$
T_1\R^3 = \{ 
(x, v) \mid x \in \R^3, v \in T_x\R^3, \| v\| = 1
\} \cong \R^3\times S^2, 
$$
to the Euclidean three space $\R^3$
has the contact structure $\{ vdx = 0 \} \subset T(T_1\R^3)$. 
We have analogous double Legendre fibrations as in the projective framework: 
$$
\begin{array}{cccccccc}
&& PT^*\R P^3 & \looparrowleft &  T_1\R^3 &   & & 
\vspace{0.5truecm}
\\
& & &  \pi_1\swarrow  &            &  \searrow\pi_2 &   & 
\vspace{0.5truecm}
\\
\R P^3 & \supset &  \R^3 &    & &
\R\times S^2 & \looparrowright & \R P^3, 
\end{array}
$$
where $\pi_1$ is the bundle projection and $\pi_2$ is defined by 
$
\pi_2(x, v) = (-x\cdot v, v), 
$
$\R\times S^2$ being identified with the space of co-oriented affine planes in $\R^3$. 
Note that $T_1\R^3$ is mapped to $PT^*(\R P^3)$ by 
$\Phi: (x, v) \mapsto ([1, x], [-x\cdot v, \ v])$ as a double covering on the image, 
that the mapping $\Phi : T_1\R^3 \to PT^*(\R P^3)$ 
is a local contactomorphism, and that $\R\times S^2$ is mapped to $\R P^3$ by 
$(r, v) \mapsto [r, v]$ as a double covering on the image which is $\R P^3 \setminus \{ [1, 0, 0, 0]\}$. 


Any co-oriented surface with boundary $(S, \gamma)$ in $\R^3$ 
lifts to a Legendre surface with boundary $(L, \Gamma)$ in $T_1\R^3$ uniquely. 
A generic surface in $\R^3$ induces a generic Legendre surface. 
The lifted Legendre surface $(L, \Gamma)$ 
projects to a front with boundary (boundary-front) in $\R\times S^2$ by $\pi_2$. 
Actually the \lq\lq local contact nature" of the double Legendre fibrations is the same, as is noted above, 
in projective and in Euclidean framework. 

\ber
{\rm 
There exists no invariant metrics on $T_1\R^3$ and on $\R\times S^2$ 
under the group $G$ of Euclidean motions on $\R^3$ compatible with the double fibration 
$\R^3 \leftarrow T_1\R^3 \rightarrow \R\times S^2$. Note that $G$ is not compact. 
In this sense, there is no dual Euclidean geometry: 
Duality in the level of Euclidean geometry is not straightforward, compared with projective geometry. 
As for related result on duality in Euclidean geometry, see \cite{BG}\cite{BF}. 
}
\enr

Let $S \subset \R^3$ be a co-oriented immersed surface with boundary $\gamma$. 

The $1$-st fundamental form $I : TS \to \R$ is defined by 
$I(v) := g_{\mbox{\tiny{\rm Eu}}}(v, v) = \Vert v\Vert^2$. 
The $2$-nd fundamental form $II : TS \to \R$ 
is defined by $II(v) := - g_{\mbox{\tiny{\rm Eu}}}(v, \nabla_v \nn)$, where 
$\nn : S \to T\R^3$ is the unit normal to $S$. 
Then we have $(I, II) : TS \to \R^2$, which determines the surface with boundary 
essentially. 

Set 
$G = {\mbox{\rm Euclid}}(\R^3) \subset {\mbox{\rm GL}}(4, \R)$, 
the group of Euclidean motions on $\R^3$. 
We consider Maurer-Cartan form of $G$,  
$$
\omega = \left(
\begin{array}{cccc}
0 & 0 & 0 & 0 \\
\omega^1 & 0 & - \omega_1^2 & - \omega_1^3
\vspace{0.2truecm}
\\
\omega^2 & \omega_1^2 & 0& - \omega_2^3
\vspace{0.2truecm}
\\
\omega^3 & \omega_1^3 & \omega_2^3 & 0 
\end{array}
\right). 
$$
For a surface with boundary, we have the adopted moving frame 
$\widetilde{\gamma} = (\gamma, \ee_1, \ee_2, \ee_3) : 
\R \to G$ by 
$\ee_1 = \gamma'$, the differentiation by arc-length parameter, $\ee_2$, the inner normal to $\gamma$, 
and $\ee_3 = \ee_1\times\ee_2 = \nn$. 
which is different from the Frenet-Serre frame. 

The structure equation is given by 
$$
d(\gamma(s), \ee_1(s), \ee_2(s), \ee_3(s)) = 
(\gamma(s), \ee_1(s), \ee_2(s), \ee_3(s))\widetilde{\gamma}^*\omega. 
$$
Thus we have 
$$
d(\ee_1, \ee_2, \ee_3) = (\ee_1, \ee_2, \ee_3)
\left(
\begin{array}{ccc}
0 & - \kappa_1 & - \kappa_2 \\
\kappa_1 & 0 & - \kappa_3 \\
\kappa_2 & \kappa_3 & 0
\end{array}
\right)
ds. 
$$
Namely we have 
$$
\left\{
\begin{array}{rcl}
\ee_1' & = & \ \  \kappa_1\ee_2 + \kappa_2\ee_3,  \\
\ee_2' & = & - \kappa_1\ee_1 + \kappa_3\ee_3,  \\
\ee_3' & = & - \kappa_2\ee_1 - \kappa_3\ee_2. 
\end{array}
\right.
$$
See \cite{IL}, for instance. 

Note that $\kappa_1 = \ee_2\cdot \gamma''$, $\kappa_2 = \ee_3\cdot \gamma''$ 
and that 
$\kappa_3 = II(\ee_1, \ee_2)$. 


\

Suppose $(S, \gamma)$ is a $C^\infty$ surface with boundary $\gamma$. 
Suppose the boundary $\gamma(t)$ is of finite type at $t = t_0$. 
Since $\gamma$ is an immersed curve, the type is written 
as $(a_1, a_2, a_3) = (1, 1+s, 1+s+r)$, for some positive integers $r, s$. 

Then we have

\bet
\label{osculating-tangent2}
Let $(S, \gamma)$ be a $C^\infty$ surface with boundary 
in $\R^3$. Suppose $\gamma$ is of finite $(1, 1+s, 1+s+r)$. 
Then $\gamma(t)$ has an osculating-tangent point at $t = t_0$ 
if and only if $\kappa_2^{(s-1)}(t_0) = 0$. 
\ent

\Proof
First remark that $\rank A_1(t) = \rank \gamma'(t) = 1$. 
Then $\rank A_2(t) = \rank (\gamma'(t), \gamma''(t)) = 1$ if and only if 
$\gamma''(t) (= \ee_1'(t))$ is a scalar multiple of $\gamma'(t) (= \ee_1)$, 
and the condition is equivalent to that $\kappa_2(t) = 0, \kappa_3(t) = 0$. 
Similarly we have that 
$\rank A_i(t) = 1, (1 \leq i \leq s)$ if and only if $\kappa_1^{(j)}(t) = 0, \kappa_2^{(j)}(t) = 0, 
(0 \leq j \leq s - 2)$. 
Then 
$$
\gamma^{(s+1)}(t) = \ee_1^{(s)}(t) = \kappa_1^{(s-1)}(t)\ee_2(t) + \kappa_2^{(s-1)}(t)\ee_3(t). 
$$
Moreover we have  $\rank\, A_{s+1}(t) = 2$ if and only if $(\kappa_1^{(s-1)}(t), \kappa_2^{(s-1)}(t)) \not= (0, 0)$. 
In this case the osculating plane is spanned by $\gamma'(t) = \ee_1(t), \gamma^{(s+1)}(t) = \ee_1^{(s)}(t)$. 
Therefore the osculating plane coincides with the tangent plane, which is spanned by $\ee_1(t), \ee_2(t)$,  
if and only if $\kappa_2^{(s-1)}(t) = 0$. 
\QED 

\

\noindent
{\it Proof of Theorem \ref{osculating-tangent}:} 
Generically $\gamma(t)$ is of type $(1, 2, 3)$ or $(1, 2, 4)$. Therefore, applying Theorem 
\ref{osculating-tangent2} in the case $s = 1$, 
the osculating-tangent point is characterised by $\kappa_2 = 0$. 
\QED

\

The flat extension problem is concerned with osculating-tangent points of the dual boundary 
$\widehat{\gamma}$, not $\gamma$. Actually we have 

\bep
\label{dual osculating-tangent}
Let $(S, \gamma)$ be a $C^\infty$ surface with boundary. 
Suppose $\widehat{\gamma}(t)$ is of type $(1, 2, 2+r)$ at $t = t_0$ for some positive integer $r$. 
Then $\widehat{\gamma}(t_0)$ is an osculating-tangent point for $(S^{\vee}, \widehat{\gamma})$ 
if and only if $\kappa_2 = 0$. Therefore, under the above condition, we have that 
$\gamma(t_0)$ is an osculating-tangent point for $(S, \gamma)$ if and only if 
$\widehat{\gamma}(t_0)$ is an osculating-tangent point for $(S^{\vee}, \widehat{\gamma})$. 
\enp

The proof of Proposition \ref{dual osculating-tangent} is given below in 
the proof of Proposition \ref{distance2}.

\

We show Theorem \ref{swallowtail-tangent} in more general context: 

\bet
\label{swallowtail-tangent2}
{\rm (The characterisation of swallowtail-tangent)} 
Let $(S, \gamma)$ be a $C^\infty$ surface with boundary of finite type in $\R^3$. 
Then we have 
\vspace{0.2truecm}
\\
{\rm (1):}  
A point on the boundary $\gamma$ 
is a swallowtail-tangent point with 
{\rm (I)} 
$\kappa_2 \not= 0$ 
if and only if 
the conditions {\rm (II), (III)} of Theorem \ref{swallowtail-tangent} hold. 
\vspace{0.2truecm}
\\
{\rm (2):} 
A point on the boundary $\gamma$ 
is a swallowtail-tangent point with $\kappa_2 = 0$ 
if and only if 
$
\ 
{\rm (I)'} 
 \ \kappa_1 \not= 0, \ \kappa_3 \not= 0,  \ 
{\rm (II)'} \ \kappa_2' = \dfrac{1}{2}\kappa_1\kappa_3, $
and 
$
{\rm (III)'} \ \kappa_3'' \not= \dfrac{4}{3}\kappa_3\kappa_1'. 
$
\ent

\ber
{\rm 
(1) A swallowtail-tangent point with $\kappa_2 = 0$ 
does not appear generically. 
\\
(2) The criteria of Theorem \ref{swallowtail-tangent2} has the similarity in the form to 
the general criterion of swallowtail found in \cite{KRUSUY}. 
}
\enr

\

\noindent
{\it Proof of Theorem \ref{swallowtail-tangent2}.} 
The dual-boundary $\widehat{\gamma}$ is given by $(\nn, -\gamma\cdot\nn) 
: (\R, 0) \to S^2\times\R$. Since $S^2\times\R$ is mapped to $\R P^{3*}\setminus \{[1,0,0,0]\}$ 
as a double covering, $\widehat{\gamma}$ is regarded as a curve in $\R P^{3*}$. 
To see the type of $\widehat{\gamma}$ 
we examine the $4\times (r+1)$ matrix
$$
\tilde{A}_r(t) = \left(
\begin{array}{ccccc}
\nn(t) & \nn'(t) & \nn''(t) & \cdots & \nn^{(r)}(t) 
\vspace{0.2truecm}
\\
-\gamma\cdot\nn(t) & (-\gamma\cdot\nn)'(t) & 
(-\gamma\cdot\nn)''(t) & 
\dots & 
(-\gamma\cdot\nn)^{(r)}(t) 
\end{array}
\right), 
$$
$(r = 1, 2, \dots).$ 
In terms of homogeneous coordinates,  
the curve $\widehat{\gamma}(t)$ is of type $(a_1, a_2, a_3)$ at $t = t_1$ 
if and only if, 
$$
\begin{array}{c}
\min \{ r \mid \rank \tilde{A}_r(t_0) = 2\} = a_1,  \quad 
\min \{ r \mid \rank \tilde{A}_r(t_0) = 3\} = a_2, 
\vspace{0.2truecm}
\\
\min \{ r \mid \rank \tilde{A}_r(t_0) = 4\} = a_3. 
\\
\end{array}
$$
In fact $\rank \tilde{A}_r(t) = \rank A_r(t) + 1$, for the matrix $A_r(t)$ introduced 
in \S \ref{Projective geometry of front-boundaries.}. 

As is mentioned, the boundary-envelope of $(S, \gamma)$, namely, 
the dual surface to the dual-boundary $\widehat{\gamma}$ has 
the cuspidal edge along the dual curve $\widehat{\gamma}^*$ of $\widehat{\gamma}$, 
where $\widehat{\gamma}(t)$ is of type $(1, 2, 3)$. 
The curve $\widehat{\gamma}$ in $\R P^{3*}$ is of type $(1, 2, 3)$ at $t = t_0$ 
if and only if $\det A_3(t_0) \not= 0$. 
In fact the condition is equivalent to that 
$$
\det (\widehat{\gamma}'(t_0), \widehat{\gamma}''(t_0), \widehat{\gamma}'''(t_0) ) \not= 0. 
$$

Similarly, the boundary-envelope of $(S, \gamma)$ 
is diffeomorphic to the swallowtail at the point $\widehat{\gamma}^*(t_1)$ 
in $\R^3 \subset \R P^3$ if and only if 
$\widehat{\gamma}(t)$ is of type $(1, 2, 4)$ at $t = t_1$. 
The condition is equivalent to that 
$$
\rank \tilde{A}_1 = 2, \ 
\rank \tilde{A}_2 = 3, \ 
\rank \tilde{A}_3 = 3, \ 
\rank \tilde{A}_4 = 4, 
$$
at $t = t_1$. 
Then, by the straightforward calculation, using the structure equation explained above, 
we have the criteria in Theorem \ref{swallowtail-tangent2}. 
In fact, from $\gamma' = \ee_1, \nn = \ee_3, \gamma'' = \kappa_1\ee_2 + \kappa_2\ee_3$, 
we have 
$$
\begin{array}{l}
\gamma'''  \, = \  - (\kappa_1^2 + \kappa_2^2)\ee_1 + (\kappa_1' - \kappa_2\kappa_3)\ee_2 
+ (\kappa_2' + \kappa_1\kappa_3)\ee_3, 
\\
\gamma''''  = 
(-3\kappa_1\kappa_1' - 3\kappa_2\kappa_2')\ee_1 
+ (- \kappa_1^3 - \kappa_1\kappa_2^2 - \kappa_1\kappa_3^2 + 2\kappa_1'\kappa_3 + \kappa_1\kappa_3' + 
\kappa_2'')\ee_3
\end{array}
$$
Moreover we have 
$\gamma''\cdot\nn' = \kappa_2, \ \gamma''\cdot\nn' = - \kappa_1\kappa_3$. 
Thus we have 
$$
\begin{array}{ccl}
(\gamma\cdot\nn)' - \gamma\cdot\nn' & = &  = 0, 
\\
(\gamma\cdot\nn)'' - \gamma\cdot\nn'' & = & 
- \kappa_2, 
\\
(\gamma\cdot\nn)''' - \gamma\cdot\nn''' & = & 
- 2\kappa_2' + \kappa_1\kappa_3
\\
(\gamma\cdot\nn)'''' - \gamma\cdot\nn'''' 
& = & 
\kappa_1^2\kappa_2 + \kappa_2^3 + \kappa_2\kappa_3^2 + 
2\kappa_3\kappa_1' + 3\kappa_1\kappa_3' - 3\kappa_3''. 
\end{array}
$$
Then the condition $\rank \tilde{A}_1(t_1) = 2$ is equivalent to that $\kappa_2 \not= 0, 
\kappa_3 \not= 0$ at $t = t_1$. 
The condition 
$\rank \tilde{A}_2(t_1) = 3$ is equivalent to 
that $\kappa_2 \not= 0$, or $\kappa_2 = 0, \kappa_3 \not= 0, 
\kappa_2\kappa_3 - \kappa_2' \not= 0$ at $t = t_1$. 

Let us see the condition $\rank \tilde{A}_3(t_1) = 3$, namely that $\det(\tilde{A}_3(t_1)) = 0$. 
We set $D = \det(\tilde{A}_3(t_1))$. Then 
we have, after simplifying the determinant and taking the transpose of $\tilde{A}_3$, 
$$
D = 
\left\vert
\begin{array}{cc}
\ee_3 & 0 \\
- \kappa_2\ee_1 - \kappa_3\ee_2 & 0 \\
(\kappa_1\kappa_3 - \kappa_2')\ee_1 + (- \kappa_1\kappa_2 - \kappa_3')\ee_2 & \kappa_2 
\\
A\ee_1 + B\ee_2 & 2\kappa_2' - \kappa_1\kappa_3
\end{array}
\right\vert, 
$$
where we set $\nn''' = A\ee_1 + B\ee_2 + C\ee_3$, 
$$
\begin{array}{rcl}
A & = & \kappa_3\kappa_1' + 2\kappa_1\kappa_3' - \kappa_2'' + 
\kappa_1^2\kappa_2 + \kappa_2^3 + \kappa_2\kappa_3^2, 
\\
B & = & -\kappa_2\kappa_1' - 2\kappa_1\kappa_2' - \kappa_3'' 
+ \kappa_1^2\kappa_3 + \kappa_2^2\kappa_3 + \kappa_3^2, 
\\
C & = & - 3\kappa_2\kappa_2' - 3\kappa_3\kappa_3'. 
\end{array}
$$
Then we see that $D$ is equal to the left hand side of the condition (II) of 
Theorem \ref{swallowtail-tangent}. 

To see the condition $\rank \tilde{A}_4(t_1) = 4$ we calculate the 
sub-determinant $E$ obtained by deleting the fourth column from $\tilde{A}_4(t_1)$. 
The condition is equivalent to $E \not= 0$. 
The sub-determinant $E$ is given by 
$$
E = - 
\left\vert
\begin{array}{ccc}
\kappa_2 & \kappa_2' - \kappa_1\kappa_3 & A' - B\kappa_1 \\
\kappa_3 & \kappa_3' + \kappa_1\kappa_2 & B' - A\kappa_1 \\
0 & \kappa_2 & 
\kappa_1^2\kappa_2 + \kappa_2^3 + \kappa_2\kappa_3^2 + 
2\kappa_3\kappa_1' + 3\kappa_1\kappa_3' - 3\kappa_3''
\end{array}
\right\vert, 
$$
and it is equal to, up to sign, the left hand side of (III). 
Thus we have (1). 

To see (2), suppose the non-generic condition $\kappa_2 = 0$. 
Then we have $\kappa_3 \not= 0, \kappa_1\kappa_3 - \kappa_2' \not= 0$. 
From the condition $D = 0$ we have 
$$
(*) \quad \kappa_1^2\kappa_3^2 - 3\kappa_1\kappa_2\kappa_2' + 2(\kappa_2')^2 = 0. 
$$
From the condition $E \not= 0$, we have 
$$
(**) \quad 2\kappa_1\kappa_3^2\kappa_1' - 3(\kappa_1\kappa_3 - \kappa_2')\kappa_3'' \not= 0. 
$$
Since the equation $(*)$ has solutions $\kappa_2' = \kappa_1\kappa_3, \dfrac{1}{2}\kappa_1\kappa_3$, 
we have $\kappa_2' = \dfrac{1}{2}\kappa_1\kappa_3$ and $\kappa_1 \not= 0$. 
Then the condition $(**)$ is equivalent to that $\kappa_3'' \not= 
\dfrac{4}{3}\kappa_3\kappa_1'$. 
\QED

\

\noindent
{\it Proof of Remark \ref{osculating-tangent-elliptic}:} 
At an osculating point on $\gamma$, the second fundamental form $II$ of $S$ satisfies 
$II(\ee_1, \ee_1) = - \ee_3'\cdot\ee_1 = \kappa_2 = 0$ and 
$II(\ee_1, \ee_2) = \kappa_3$. 
Therefore $\det(II) = - \kappa_3^2 \leq 0$. 
Moreover $\det(II) < 0$ if and only if $\kappa_3 \not= 0$. 
\QED

The envelope-swallowtail point for $(S, \gamma)$ 
corresponds to the osculating plane to $\widehat{\gamma}$ at a point $t = t_1$ of type $(1, 2, 4)$ 
in $\R P^{3*}$. 
Then $d$ is the distance between $\gamma(t_1)$ and $\widehat{\gamma}^*(t_1)$. 
Actually the formula in Proposition \ref{distance} gives 
${\mbox{\rm{dist}}}(\gamma(t), \widehat{\gamma}^*(t))$: 

\bep
\label{distance2}
Let $(S, \gamma)$ be a $C^\infty$ surface with boundary. 
Suppose the dual-boundary $\widetilde{\gamma}(t)$ is of type $(1, 2, 2+r)$ 
at $t = t_1$ for some positive integer $r$. 
Then the distance $d = 
{\mbox{\rm{dist}}}(\gamma(t_1), \widehat{\gamma}^*(t_1))$ 
is given by 
$$
d = 
\left\vert 
\dfrac{\kappa_2\sqrt{\kappa_2^2 + \kappa_3^2}}{\kappa_2(\kappa_3' + \kappa_1\kappa_2) 
+ \kappa_3(-\kappa_2' + \kappa_1\kappa_3)}
\right\vert. 
$$
at $t = t_1$. 
\enp

Proposition \ref{distance} follows from Proposition \ref{distance2}. 

\

\noindent
{\it Proof of Proposition \ref{dual osculating-tangent} and Proposition \ref{distance2}:} 
Let $\gamma(t)$ be a point on the boundary $\gamma$ in $\R P^3$.  
Set $\widehat{\gamma}(t) = (- \gamma(t)\cdot\nn(t), \nn(t))$, where $t$ is the arc-length 
parameter.  
Since $\widehat{\gamma}$ is of type $(1, 2, 2+r)$, $\widehat{\gamma}', \widehat{\gamma}''$ 
are linearly independent. 
Then the point $\widehat{\gamma}^*(t) = [1, x] = [1, x_1, x_2, x_3] \in \R^3 \subset \R P^3$
is obtained 
by solving the system of equations
$$
\left\{ 
\begin{array}{ccc}
\xxx\cdot\nn - \gamma(t)\cdot\nn(t) & = & 0, \\
\xxx\cdot\nn' - (\gamma(t)\cdot\nn(t))' & = & 0, \\
\xxx\cdot\nn'' - (\gamma(t)\cdot\nn(t))'' & = & 0. 
\end{array}
\right.
$$
We set $\Delta = \left\vert
\begin{array}{ccc}
n_1  & n_2  & n_3 \\
n_1'  & n_2'  & n_3' \\
n_1''  & n_2''  & n_3''
\end{array}
\right\vert
$, 
where we set $\nn(t) = (n_1(t), n_2(t), n_3(t))$. 
Note that, under the assumption, 
the Gauss mapping $\nn(t)$ 
restricted at $\gamma$ is immersive and therefor $\Delta \not= 0$. 
Then, by Cram\'{e}r's formula, we have 
$$
\begin{array}{c}
x_1 = \dfrac{1}{\Delta}
\left\vert
\begin{array}{ccc}
\gamma\cdot\nn & n_2 & n_3 \\
(\gamma\cdot\nn)' & n_2' & n_3' \\
(\gamma\cdot\nn)'' & n_2'' & n_3''
\end{array}
\right\vert
, 
x_2 = \dfrac{1}{\Delta}
\left\vert
\begin{array}{ccc}
n_1 & \gamma\cdot\nn & n_3 \\
n_1' & (\gamma\cdot\nn)'& n_3' \\
n_1'' & (\gamma\cdot\nn)''& n_3''
\end{array}
\right\vert
, 
\vspace{0.2truecm}
\\
x_3 = \dfrac{1}{\Delta}
\left\vert
\begin{array}{ccc}
n_1  & n_2  & \gamma\cdot\nn \\
n_1'  & n_2'  & (\gamma\cdot\nn)' \\
n_1''  & n_2''  & (\gamma\cdot\nn)''
\end{array}
\right\vert
. 
\end{array}
$$
Since 
$$
(\gamma\cdot\nn)' = \gamma'\cdot\nn + \gamma\cdot\nn' = \gamma\cdot\nn', 
\ \  
(\gamma\cdot\nn)'' = \gamma'\cdot\nn' + \gamma\cdot\nn'', \ \ 
\gamma'\cdot\nn' = - \kappa_2, 
$$
we have 
$$
x_1 = \gamma_1 + 
\dfrac{1}{\Delta}
\left\vert
\begin{array}{ccc}
0 & n_2 & n_3 \\
0 & n_2' & n_3' \\
\gamma'\cdot\nn' & n_2'' & n_3''
\end{array}
\right\vert
= \gamma_1 - \dfrac{\kappa_2}{\Delta}
\left\vert
\begin{array}{cc}
n_2 & n_3 \\
n_2' & n_3' 
\end{array}
\right\vert. 
$$
Similarly we have 
$$
x_2 = \gamma_2 
+ \dfrac{\kappa_2}{\Delta}
\left\vert
\begin{array}{cc}
n_1 & n_3 \\
n_1' & n_3' 
\end{array}
\right\vert
, 
\quad 
x_3 = \gamma_3 
- \dfrac{\kappa_2}{\Delta}
\left\vert
\begin{array}{cc}
n_1 & n_2 \\
n_1' & n_2' 
\end{array}
\right\vert. 
$$
The distance $d$ between a point $\gamma(t)$ on the boundary 
and the point $\widehat{\gamma}^*(t)$ on the boundary-envelopes is calculated by 
$$
\begin{array}{rcl}
d^2  & = & 
(x_1 - \gamma_1)^2 + (x_2 - \gamma_2)^2 + (x_3 - \gamma_3)^2 
\vspace{0.2truecm}
\\
 & = & \dfrac{\kappa_2^2}{\Delta^2}\left(
\left\vert
\begin{array}{cc}
n_2 & n_3 \\
n_2' & n_3' 
\end{array}
\right\vert^2
+
\left\vert
\begin{array}{cc}
n_1 & n_3 \\
n_1' & n_3' 
\end{array}
\right\vert^2
+
\left\vert
\begin{array}{cc}
n_1 & n_2 \\
n_1' & n_2' 
\end{array}
\right\vert^2
\right). 
\end{array} 
$$
Now, from the structure equation, we have 
$$
\Delta = \vert \ee_3, \ee_3', \ee_3''\vert = \kappa_2(\kappa_3' + \kappa_1\kappa_2) + 
\kappa_3(- \kappa_2' + \kappa_1\kappa_3). 
$$
On the other hand, for the exterior product, 
$$
\nn\times\nn' (= *(\nn\wedge\nn')) 
= \ee_3\times(-\kappa_2\ee_1 - \kappa_3\ee_2) = \kappa_3\ee_1 - \kappa_2\ee_2, 
$$
$$
\vert \nn\times\nn' \vert^2 = \vert \kappa_3\ee_1 - \kappa_2\ee_2\vert^2 = 
\kappa_2^2 + \kappa_3^2. 
$$
Therefore 
$$
d^2 = \dfrac{\kappa_2^2(\kappa_2^2 + \kappa_3^2)}{[
\kappa_2(\kappa_3' + \kappa_1\kappa_2) + 
\kappa_3(- \kappa_2' + \kappa_1\kappa_3)]^2}. 
$$
Hence we have the formula of Proposition \ref{distance2} and therefore Proposition \ref{distance}. 
Moreover, we see $\widehat{\gamma}^*(t_1)$ coincides with $\gamma(t_1)$ 
if and only if $\kappa_2^2(\kappa_2^2 + \kappa_3^2) = 0$, which 
is equivalent to that $\kappa_2 = 0$ at $t = t_1$. 
Thus we have Proposition \ref{dual osculating-tangent}. 
\QED

To show Theorem \ref{solution of problem}, 
we show first

\bel
\label{uniqueness}
Let $(S, \gamma)$ be a $C^\infty$ surface with boundary, 
$\widetilde{S}$ a flat $C^1$ extension of $S$. 
Suppose the restriction $g\vert \gamma$ of the Gauss mapping of $S$ restricted on $\gamma$ 
is an immersion. 
Then for any $p \in \gamma$, there is an open neighbourhood $U$ of $p$ in 
$\widetilde{S} \setminus {\mbox{\rm{Int}}}S$ such that 
the Legendre lifting of $U$ 
projects by $\pi_2$ to $\widehat{\gamma}$. 
\enl

\Proof
Set $S' = \widetilde{S} \setminus {\mbox{\rm{Int}}}S$. 
Note that $S'$ is a $C^\infty$ surface with boundary $\gamma$. 
Consider the Legendre liftings $(L, \Gamma)$ of $(S, \gamma)$ and 
$(L', \Gamma)$ of $S'$ in the incident manifold 
$I^5$ with respect to the projection $\pi_1$. 
Because $\widetilde{S}$ is a $C^1$ surface, 
we see $\widetilde{L} = L \cup L'$ is a $C^0$ surface in $I$. 
From the assumption that $g\vert_{\gamma}$ 
is immersive, we see the $S^2$ component of $\pi_2\vert_{\Gamma} : \Gamma \to S^2\times \R$ 
is immersive. 
Consider the Gauss mapping $g'$ of $S'$ and its restriction $g'\vert_{\gamma}$. 
Then $g'\vert_{\gamma}$ is immersive if and only if the $S^2$-component 
of $\pi_2\vert_{\Gamma}$ is immersive. Since $S'$ is flat, $g'$ is of rank $< 2$. 
Hence $g'$ is of rank one along $\gamma$. Therefore $\pi_2\vert_{L'}$ 
is of rank one along $\Gamma$. Moreover the kernel field of $\pi_2\vert_{L'}$ 
is transverse to $\Gamma$ on $\Gamma$. 
Then $L'$ projects to $\widehat{\gamma}$ near $\Gamma$. 
\QED

\

\noindent
{\it Proof of Theorem \ref{solution of problem}:} 
Suppose $(S, \gamma)$ is generic. Then 
the dual-boundary $\widehat{\gamma}$ is of type $(1, 2, 3)$ or $(1, 2, 4)$ 
(Theorem \ref{basic results}). 
Suppose $p \in \gamma$ is not an osculating-tangent point. 
the boundary-envelope $E$ is non-singular near $p$ (Theorem \ref{osculating-tangent2}). 
Then, actually, the pair $(S, E)$ is of type $B_2$ and $(S, \gamma)$ has the $C^1$ flat 
extension by $E$.  
To show the uniqueness of local flat extensions, 
suppose $(S, \gamma)$ has a local $C^1$ flat extension $\widetilde{S}$. 
Then by Lemma \ref{uniqueness}
the Legendre lifting $L'$ of $\widetilde{S} \setminus {\mbox{\rm{Int}}}S$ 
projects to $\widehat{\gamma}$ locally at each point of $\Gamma$. 
Therefore $L'$ is contained in the projective conormal bundle 
of $\pi_2\vert_{\Gamma}$. Hence, by projecting by $\pi_1$, 
we see that $\widetilde{S} \setminus {\mbox{\rm{Int}}}S$ is locally contained in 
the boundary-envelope $\pi_1(L_2)$. Thus we have the local uniqueness 
of the flat extension. 
\QED 

{

}

\

\

\noindent Goo ISHIKAWA \\
Department of Mathematics, \\
Hokkaido University,\\
Sapporo 060-0810, JAPAN. \\
{\vspace{-0.7truecm}}
\begin{verbatim}
E-mail : ishikawa@math.sci.hokudai.ac.jp
\end{verbatim}

\end{document}